\documentclass[intlimits,a4paper]{amsart}
\usepackage[all]{xy}
\usepackage{amssymb}

\newcommand{\R}{{\bf R}}

\newcommand{\Q}{{\bf Q}}

\newcommand{\Z}{\bold{Z}}
\newcommand{\N}{\bold{N}}
\newcommand{\C}{\bold{C}}

\newcommand{\tr}{\text{tr}}

\newcommand{\Sym}{\text{Sym}}

\newtheorem{theorem}{Theorem}[section]

\newtheorem{corollary}[theorem]{Corollary}
\newtheorem{proposition}[theorem]{Proposition}
\numberwithin{equation}{section}
\newcommand{\abs}[1]{\lvert#1\rvert}

\begin{document}
\title[Local theta correspondence and the DII lifting]{Local theta correspondence and the lifting of Duke, Imamoglu
  and Ikeda}  
\author[R. Schulze-Pillot]{Rainer Schulze-Pillot} 
\thanks{MSC 2000: Primary 11F46, Secondary 11F27}
\maketitle
\begin{abstract}
{We use results on the local theta correspondence to prove that for
  large degrees the Duke-Imamoglu-Ikeda lifting of an elliptic modular
form is not a linear combination of theta series.}
\end{abstract}

\section{Introduction} 
In the article \cite{ksm} it is proved as a side remark that a Siegel
modular form $F$ of degree $2n$ and weight $n+k$ obtained by the lifting
of Duke, Imamoglu and Ikeda 
\cite{ikeda} (called the DII-lifting in the sequel) from
an elliptic modular form of weight $2k$ is not
a linear combination of theta series of even unimodular positive
definite quadratic forms of rank $m=2(n+k)$ if  $n$ is bigger than
$k$, whereas for $n=k\equiv 0 \bmod 2$   
the DII liftings lie in the space generated by theta series subject to
a conjecture on $L$-functions of elliptic cuspidal Hecke
eigenforms. The proof uses B\"ocherer's characterization of the cuspidal
Siegel eigenforms that lie in the space of theta series by special
values of their standard $L$-functions.

In this article we use results on the local theta correspondence to
give a different proof of the first result. In fact we prove  a more general
version including theta
series attached to arbitrary non degenerate quadratic forms and also
theta series with spherical harmonics. The underlying local fact has
been noticed by several people (including the author) immediately
after the preprint 
version of \cite{ikeda} became available in 1999 but has apparently never been
published. 

In the case
$n=k$ we show that the local representations attached
to a DII-lift 
are in the image of the local theta correspondence with the orthogonal
group of  a suitable  local quadratic space of dimension $4n=2(n+k)$
for all places. 

I thank Ralf Schmidt for discovering a mistake in an earlier version
of this article.
 
\section{The case $n>k$}
Let $f$ be an elliptic modular form of weight $2k$ for the full
modular group $SL_2(\Z)$. It was conjectured by Duke and Imamoglu and
proven by Ikeda in \cite{ikeda} that for any $ n\equiv k \bmod 2$
there exists a nonzero Siegel cusp form $F=F_{2n}(f)$ of weight $n+k$ for
the group $Sp_{2n}(\Z) \subseteq SL_{4n}(\Z)$ whose standard
$L$-function is equal to
$$\zeta(s)\prod_{i=1}^{2n}L(s+k+n-i,f),$$
where $L(s,f)$ is the usual Hecke $L$-function of $f$. We call
$F_{2n}(f)$ the DII-lift of degree $2n$ of $f$.

In fact, it has been shown in \cite[Lemma 1.3.1]{schmidt},
that the Satake
parameters $\alpha_p^{(j)}$ of $F$ at the prime $p$ are given by 
\begin{equation}\label{satakepar1}
\alpha_p^{(0)}=\beta_p^{-n}, \quad\alpha_p^{(j)}=\beta_p
p^{-n+j-\frac{1}{2}} \quad\text{for } 1\le j \le 2n,
\end{equation}
where we write the 
$p$-factor of $L(s,f)$ as $(1-\beta_p
p^{-s+k-\frac{1}{2}})^{-1}(1-\beta_p^{-1} p^{-s+k-\frac{1}{2}})^{-1}$,
i. e., we have (by Ramanujan-Petersson) 
the normalization $\abs{\beta_p}=1$.

Moreover, Ikeda has recently announced a generalization of this
result. In this generalized version $f$ is replaced by an irreducible cuspidal
automorphic representation $\tau=\otimes \tau_v$ of $GL_2({\bf A}_E)$,
where $E$ is a 
totally real number field and $\tau_v$ is assumed to be a principal
series representation attached to a character $\mu_v$ for
all finite places $v$ of $E$ and to be 
discrete series with minimal weight $\pm \kappa_w$ for the infinite
places $w$ of $E$. 

Ikeda then proves the existence of an irreducible
cuspidal automorphic representation $\pi=\pi(m,\tau)$ of the
metaplectic group $\widetilde{Sp_{m}}({\bf A}_E)$ with local components $\pi_v$
described in terms of $\tau_v$ as follows: For a real place $w$ the 
representation $\pi_w$ is the lowest weight representation of lowest
$K$-type $\det^{\kappa_w+m/2}$, for a finite place $v$ the 
representation $\pi_v$ is a degenerate principal series representation
which is induced from a character $\mu_v^{(m)}$ on the
maximal parabolic $\widetilde{P_m}$ derived from $\mu_v$

If $m=2n$  is even one can obtain from $\pi$ a representation of
$Sp_{2n}({\bf A}_E)$, also denoted by $\pi=\pi(2n, \tau)$; this
representation is induced from the $2n$-tuple of characters 
$\mu_v \abs{\quad}^{-n+j-\frac{1}{2}}$ for $1\le j \le 2n$. 
If (for
$E={\bf Q}$) in addition $f$ is as above 
and $\tau$ is the automorphic representation of $GL_2({\bf
  A}_{\Q})$ associated to $f$ this gives the representation of
$Sp_{2n}({\bf   A}_{\Q})$ associated to $F_{2n}(f)$. 

\vspace{0.5cm}
We are interested in the question whether the Siegel modular form
$F_{2n}(f)$ can be
obtained as a linear combination of theta series, respectively in the
more general situation whether the representation $\pi(2n,\tau)$ is
in the image of the theta correspondence between $Sp_{2n}$ and a
suitable orthogonal group. 

We let $V$ be a vector space over $E$ of even dimension $2r$ with
a nondegenerate quadratic form $q$ and associated symmetric bilinear
form $B(x,y)=q(x+y)-q(x)-q(y)$ on it. If $E={\bf Q}$  and $q$
is positive definite we consider a homogenous pluriharmonic form $P\in
\C[\{X_{ij}\mid 1 \le i\le 2r, 1\le j \le m\}]$ in $2rm$ variables of
weight  $\nu$ (see \cite[III, 3.5]{Freitag}) as a function on
$(V\otimes \R)^m$ by identifying the 
latter space with $\R^{2rm}$ using a basis of $V$ which is orthonormal
with respect to the symmetric bilinear form $B$.

The theta series of degree $m$  of $L$ with respect to $P$
is defined for any   ${\bf Z}$ - lattice $L$ of full rank $2r$
on $V$  by 
$$\vartheta^{(m)}(Z,L,P)=\sum_{{\bf x}\in 
L^m}P({\bf x})\exp(2\pi i \tr(q({\bf x})Z)),$$ 
where $Z\in {\mathfrak H}_m$ is in
the Siegel upper half space ${\mathfrak H}_{m}\subseteq M_{m}^{\Sym}(\C)$
and where $q({\bf x})=\frac{1}{2}(B(x_i,x_j))$ is half the Gram matrix
of the $m$-tuple ${\bf x}=(x_1,\ldots,x_{m})$.  It is a Siegel modular
form of weight $r+\nu$ for some congruence subgroup of $Sp_m(\Z)$. One
has a similar 
definition for arbitrary totally real $E$ and totally positive
definite $q$.

If $L$ is even unimodular the theta series is a Siegel modular form
for the full modular group $Sp_{m}(\Z)$.
 
More generally one considers for arbitrary $E$ and (non degenerate)
$q$ the theta correspondence associating to a subset of the set of
irreducible automorphic representations of the adelic
orthogonal group  $O_{(V,q)}({\bf A}_E)$ a subset of the set of
irreducible automorphic representations of the adelic group 
$Sp_m({\bf A}_E)$. There is by now a vast literature (starting with
Weil's \cite{Weil} and Howe's \cite{Howe}) but no textbook
reference on this correspondence, for definitions and properties of it
see \cite{rallis_functoriality,moeglin1,moeglin2,kudla_local,mvw}.

\newpage
\begin{theorem}
Let $E,V,q,\tau$ be as above.
\begin{enumerate}
\item If $2n>r-1$ and there is a finite place $v$ of $E$ for which the completion
  $V_v$ of the quadratic space 
$(V,q)$ is not split (i. e. is not an orthogonal sum of hyperbolic planes)
the representation $\pi(2n, \tau)$ is not in the image of the theta correspondence with
$O_{(V,q)}({\bf A}_E)$. 
\item 
If $2n>r$  the representation
$\pi(2n,\tau)$ is not in the image of the theta correspondence with
$O_{(V,q)}({\bf A}_E)$. 
\end{enumerate}
\end{theorem}
\begin{corollary}
Let $f$ be an elliptic modular form of weight $2k$ as above, $\nu \in
\N_0$.\\ 
Then for $n>k-\nu$ the DII-lift $F_{2n}(f)$  is not a linear combination
of theta series of positive definite quadratic forms with pluriharmonic
forms of degrees $\nu'\ge \nu$.

In particular for $n>k$ the DII-lift $F_{2n}(f)$ is not a linear
combination of theta 
series attached to positive definite quadratic forms (with or without
pluriharmonic forms).
\end{corollary}
\vspace{0.3cm}
\begin{proof}[Proof of the Theorem]
If $\pi(2n,\tau)$ is in the image of the global theta correspondence
with $O_{(V,q)}({\bf A}_E)$ its local components $\pi_v(2n,\tau)$ are in the
image of the local theta correspondence with  $O_{(V,q)}(E_v)$ for all
places $v$ of $E$; we denote by $\pi_v'$ the corresponding representation
of this orthogonal group.  

The local theta correspondence 
has been described explicitly in terms of the Bern\-stein--Zelevinsky
data of the representations in
\cite{rallis_functoriality,kudla_local}.
By construction the Bern\-stein--Zelevinsky data of $\pi_v(2n,\tau)$ are
$\mu_v \abs{\quad}^{-n+j-\frac{1}{2}}$ for $1\le j \le 2n$. 
The characters $\mu_v
\abs{\quad}^{-n+j-\frac{1}{2}}$ are never of the type $\abs{\quad}^i$
for some integral $i$; this is clear if $\mu_v$ is ramified, and
follows from $\abs{\mu_v(\omega_v)}<q_v^{\frac{5}{34}}$ (see
\cite{kimshahidi}) if $\mu_v$ is unramified (where $\omega_v$ is a
prime element at the place $v$ and $q_v$ is the norm of $\omega_v$).
 Since by the results of \cite{rallis_functoriality,kudla_local} only
 characters of the type  $\abs{\quad}^i$
for some integral $i$ can be missing in the 
Bernstein--Zelevinsky data of the
 representation $\pi'$ on the orthogonal side
 we see that
 all the $\mu_v \abs{\quad}^{-n+j-\frac{1}{2}}$ for $1\le j
 \le 2n$ have to appear in these Bernstein--Zelevinsky data.
Obviously the rank
 of the orthogonal group  $O_{(V,q)}(E_v)$ has to be at least $2n$ for
 this to be possible. If $V_v$ is split this requires $r\ge 2n$, if
 $V_v$ is not split it requires $r-1\ge 2n$ (and even $r-2 \ge 2n$ if
 the anisotropic kernel of $V_v$ is of dimension $4$).

\end{proof}
\begin{proof}[Proof of the Corollary]
It is well known that one can associate to a cuspidal Siegel modular
form $F$ for
$Sp_{2n}(\Z)$ which is an eigenfunction of all Hecke operators an
irreducible cuspidal automorphic representation $\pi(F)$ of
$Sp_{2n}({\bf   A}_\Q)$ with the same Satake parameters, see
e.g. \cite{Asgari-Schmidt}. 
Moreover, if $F$ has weight $\tilde{k}$ and can be written as a
linear combination of theta series of positive definite quadratic
forms with pluriharmonic forms of weights $\nu'\ge \nu$, 
the fact that the space generated by theta series with pluriharmonic
forms of weight $\nu'$ for lattices on a fixed quadratic space is
Hecke invariant guarantees that
there is a 
positive definite quadratic space $(V,q)$ of dimension
$2r=2(\tilde{k}-\nu')$
for some $\nu'\ge \nu$
such that the irreducible representation $\pi(F)$ is in the image of
the theta correspondence with  
 $O_{(V,q)}({\bf A}_{\Q})$. 

Considering this for the  DII lift $F_{2n}(f)$ and $\nu, \nu',k$ as in the
assertion we have $\tilde{k}=k+n$ and
therefore $r=k+n-\nu'\le k+n-\nu$ and
 $k<n+\nu$, hence $r<2n$.

By the Theorem $\pi(F_{2n}(f))$ is not in
 the image of the theta correspondence with  $O_{(V,q)}({\bf
   A}_{\Q})$, so $F_{2n}(f)$ can not be a linear combination of theta
 series as described above.
\end{proof}

\newpage
{\it Remark.}
\begin{enumerate} \item If we omit the restriction that the quadratic
  form is positive 
  definite we still get the same result as long as we restrict
  attention to the holomorphic theta series of weight $r$ associated
  to indefinite quadratic forms of rank $2r$ that have
  been constructed by Siegel and Maa\ss\ \cite{siegel,maass}. Although
  there is no method known  to use suitable test functions
  in the oscillator representation in order to construct
holomorphic theta series of weight $k+n-\nu'$ associated (by the theta
correspondence) to
indefinite quadratic forms of rank larger than  $2(k+n-\nu')$, the
results of Rallis \cite{rallis_howeduality} seem 
to imply that such a construction is indeed possible.
Such a 
construction would then yield $F_{2n}(f)$ without contradicting the
Theorem.
\item The Corollary could in principle be proved without using
  representation theoretic tools by computing the possible Hecke
  eigenvalues of a linear combination of theta series with the help of
  results of Andrianov \cite{andrianov_hecke} and Yoshida
  \cite{yoshida} and comparing with the eigenvalues of a DII-lift; we
  expect this to be a rather tedious and unpleasant computation.
\end{enumerate}
\section{The case $n=k$}
With the notations of the previous section we assume now $E=\Q$ and
$n=k$, so that the weight $k+n$ of a DII lift is equal
to the rank $2n$ of the symplectic group considered.
\begin{theorem}
For $n=k$ the local components $\pi_p$ of the representation $\pi(F_{2n}(f))$ of
the Duke-Imamoglu-Ikeda lift $F_{2n}(f)$ are in the image of the local theta
correspondence with the split quadratic space over $\Q_p$ of dimension
$4n=2(k+n)$ (which is the orthogonal sum of $2n$ hyperbolic planes) for all
(finite) primes $p$.

The component $\pi_\infty$ at the real place is in the image of the
theta correspondence with the orthogonal group of the positive definite quadratic space over
$\R$ of dimension $4n$ and also  in the image of the
theta correspondence with the orthogonal group of the quadratic space
of dimension $4n$ and signature $(4n-1,1)$ over $\R$.
\end{theorem}
\begin{proof}
The assertion for the finite primes is again an immediate consequence
of the results of \cite{rallis_functoriality,kudla_local}.
The assertion at the real place follows from Theorem 15 of \cite{paul}
since $\pi_\infty$ is a limit of discrete series.
\end{proof}
{\it Remark.}
\begin{enumerate}
\item At the real place it follows from Theorem 15 of \cite{paul} that
  $\pi_\infty$ occurs also in the images of the
theta correspondence with the orthogonal groups of the quadratic spaces
of dimension $4n+2$ and signatures $(4n+1,1)$ and $(4n,2)$ over $\R$,
but not for any other quadratic space of dimension $4n$ or $4n+2$.
\item Since the split quadratic space of
  dimension $4n$ has square discriminant, there is  no
  quadratic space over $\Q$ which is split at all finite primes and of
  signature   $(4n-1,1)$. It is moreover well known that there is a
  positive definite quadratic space over $\Q$ of even dimension $2r$ which is split at all
  finite primes if and only if $r$ is divisible by $4$ and that the
  same condition is necessary and sufficient for the existence of a
  quadratic space over $\Q$ of even dimension $2r+2$ and signature
  $(2r+1,1)$ which is split at all
  finite primes

\parindent=0pt
We see that $\pi(F_{2n}(f))$ can (in the case $n=k$) be in the image of
the global theta correspondence with some quadratic space of dimension
$4n=2(n+k)$ only if $n$ is even and the quadratic space
is the unique positive definite space of that dimension which is split
at all finite places. 

\parindent=0pt
In a similar way for even $n$ there is also a unique quadratic
space $(V,q)$ over $\Q$ of dimension 
$4n+2$ and signature $(4n+1,1)$ for which $\pi_v$ is in the
image of the local  theta correspondence with the orthogonal group of
the completion $V_v$ for all (finite or infinite) places $v$ of $\Q$;
this space is again split at all finite primes.

There are also spaces of signature $(4n,2)$ for which $\pi_v$ is in the
image of the local  theta correspondence with the orthogonal group of
the completion $V_v$ for all (finite or infinite) places $v$ of $\Q$;
such a space can be obtained as the orthogonal sum of an imaginary
quadratic field 
equipped with the norm form scaled by some negative number and a
positive definite space of 
dimension $4n$ which is split at all finite places.
\parindent=0pt
In all cases our purely local methods allow no statement about
occurrence in the global theta correspondence for the respective space.
\end{enumerate}

\section{Iterated theta liftings}
One might ask whether it is possible to construct the
Duke--Imamoglu--Ikeda lift 
  $F_{2n}(f)$ or the associated automorphic representation by a sequence
  of theta liftings between groups $G_i$, where for each $i$ the pair 
$G_i,G_{i+1}$ consists (in either order) of a symplectic or
metaplectic group and an orthogonal group, starting with the
representation associated to $f$ on $SL_2$ 
or the representation on the metaplectic group $\widetilde{SL_2}$
associated to the form $g$ which corresponds to $f$ under the Shimura
correspondence. We have more generally:
\begin{proposition}
Let the notations be as in Sections 1 and 2 and $n>1$. 

The generalized Duke--Imamoglu--Ikeda lift $\pi(2n,\tau)$ can not be
constructed by a series of theta liftings as described above.
\end{proposition}
\begin{proof}
If one starts out with the metaplectic group $\widetilde{SL_2}$ all
subsequent groups $G_i$ will be orthogonal groups of quadratic spaces
in odd dimension or metaplectic groups  $\widetilde{Sp_m}$ (with a
genuine representation of $\widetilde{Sp_m}$ on it) so that we
will never arrive at $Sp_{2n}$, with the exception that the initial
step may lead to $SO(3,2)$ identified with $PGSp_2$ or to $SO(2,1)$
identified with $PGL_2$ (in which case we obtain a Saito--Kurokawa
lift respectively a Shimura lift). After this initial step there have
to appear correspondences which raise the rank of the group, and again by
\cite{rallis_functoriality,kudla_local} all 
local representations occurring will (at each finite place $v$ of $E$)
have at most two terms
$\mu_v\abs{\quad}^j$ with $j\in \frac{1}{2}\Z \setminus \Z$ among
their Bernstein--Zelevinsky data, so $\pi_v(2n,\tau)$ can never occur
for $n>1$.  

If one starts with $SL_2$ the representation $\pi_v(2n,\tau)$ can
not occur if $n>1$  for the same reason. 
\end{proof}

Rainer Schulze-Pillot\\
Fachrichtung 6.1 Mathematik,
Universit\"at des Saarlandes (Geb. E2.4)\\
Postfach 151150, 66041 Saarbr\"ucken, Germany\\
email: schulzep@math.uni-sb.de


\begin{thebibliography}{MVW}
\bibitem{andrianov_hecke} A. N. Andrianov: Quadratic forms and Hecke
  operators. 
{\it Grundlehren der Mathematischen Wissenschaften}, {\bf
  286}. Springer-Verlag, Berlin, 
1987.
\bibitem{Asgari-Schmidt} M. Asgari, R. Schmidt: {\it Siegel modular forms and
  representations}, Manuscripta Math. {\bf 104} (2001), 173--200. 
\bibitem{Freitag} E. Freitag: Siegelsche Modulfunktionen,  {\it Grundlehren
der Mathematischen Wissenschaften}, {\bf 254}, Springer-Verlag, Berlin,
1983. 
\bibitem{Howe} R. Howe: {\it $\theta $-series and invariant theory.
  Automorphic forms, representations and $L$-functions}, in:
  Proc. Sympos. Pure Math., Oregon State Univ., Corvallis, Ore.,
  1977, Part 1,  pp. 275--285, Proc. Sympos. Pure Math., XXXIII,
  Amer. Math. Soc., Providence, R.I., 1979. 
\bibitem{ikeda} T. Ikeda: {\it On the lifting of elliptic cusp forms to
  Siegel cusp forms of degree $2n$},  Ann. of Math. (2)  {\bf 154}  (2001),
  641--681 
\bibitem{kimshahidi} H. Kim, F. Shahidi: Functorial products
  for $\rm GL\sb 2{\times} GL\sb 3$ and functorial symmetric cube for
  $\rm GL\sb 2$.  C. R. Acad. Sci. Paris S\'er. I Math.  331  (2000),
  no. 8, 599--604. 
\bibitem{ksm} W. Kohnen, R. Salvati Manni: {\it Linear relations between
  theta series},   Osaka J. Math.  {\bf 41}  (2004), 353--356. 
\bibitem{kudla_local} S. Kudla: {\it On the local theta
    correspondence}, Invent. math. {\bf 83} (1986), 229-255
\bibitem{maass} H. Maa\ss: {\it Modulformen zu indefiniten quadratischen 
Formen}
Math. Scand.
{\bf 17}
(1965), 41--55. 
\bibitem{mvw} C. M\oe glin,M.-F. Vign\'eras, J.-L. Waldspurger:
  Correspondances de Howe sur un corps $p$-adique,
Lecture Notes in Mathematics, 1291. Springer-Verlag, Berlin, 1987.
\bibitem{moeglin1} C. M\oe glin:  {\it Non nullit\'e de certains rel\^evements
  par s\'eries th\'eta}, J. Lie Theory  {\bf 7}  (1997),  no. 2, 201--229.  
\bibitem{moeglin2} C. M\oe glin: {\it Quelques propri\'et\'es de base des
  s\'eries th\'eta},
  J. Lie Theory  {\bf 7}  (1997),  no. 2, 231--238 
\bibitem{paul} A. Paul: {\it On the Howe correspondence for
  symplectic-orthogonal dual pairs},  J. Funct. Anal.  {\bf 228}
(2005), 270--310  
\bibitem{rallis_functoriality}
S. Rallis: {\it Langlands functoriality and
the Weil representation}, Amer. J. of Math. {\bf 104} (1982), 469-515
\bibitem{rallis_howeduality} S. Rallis: {\it On the Howe duality
    conjecture}, Compositio. Math. {\bf  51} (1984), 
333--399
\bibitem{schmidt} R. Schmidt: {\it On the spin $L$-function of Ikeda's
  lifts},  Comment. Math. Univ. St. Pauli  {\bf 52}  (2003),  1--46.
\bibitem{siegel} C. L. Siegel: {\it Indefinite quadratische Formen 
und Modulfunktionen}, in: Studies and Essays Presented to R. Courant
on his 60th Birthday, January 8, 1948,  pp. 395--406, Interscience
Publisher, Inc., New York, 1948 (Ges. Abh. 3, 85-91). 
\bibitem{Weil} A. Weil: {\it Sur la formule de Siegel dans la th\'eorie des
  groupes classiques}, Acta Math.  {\bf 113}  (1965) 1--87. 
\bibitem{yoshida} H. Yoshida: {\it The action of Hecke operators on 
theta series},
Algebraic and topological theories (Kinosaki, 1984), 197--238, 
Kinokuniya, Tokyo, 1986.
\end{thebibliography}
\end{document}